\documentclass[reqno]{amsart}

\usepackage[english]{babel}
\usepackage{graphicx}
\usepackage[hidelinks]{hyperref}
\usepackage{amssymb}
\usepackage{amsmath}
\usepackage{fixmath}
\usepackage[foot]{amsaddr}
\usepackage{amsfonts}
\usepackage{bbm}
\usepackage{xcolor}
\usepackage[top=1in, bottom=1.25in, left=1.25in, right=1.25in]{geometry}
\usepackage{bm}
\usepackage{float}
\usepackage{subcaption}
\usepackage{overpic}
\usepackage{url}
\usepackage{comment}
\graphicspath{
	{images/}
}

\def\u{{\bm u}}
\def\n{{\bm n}}

\begin{document}
	
	\title{A non-intrusive data-driven ROM framework for hemodynamics problems}
	\author{
		M. Girfoglio $^1$         \and
		L. Scandurra $^1$         \and
		F. Ballarin $^1$}
	\author{
		G. Infantino $^2$         \and
		F. Nicolò $^3$            \and
		A. Montalto $^3$          \and
		G. Rozza $^1$             \and
		R. Scrofani $^4$          \and
		M. Comisso $^3$           \and
		F. Musumeci $^{3}$
	}
	
	\address{$^1$ SISSA, Scuola Internazionale Superiore di Studi Avanzati, Area di Matematica, mathLab Trieste, Italy.}
	\address{$^2$ Politecnico di Torino, Collegio di Ingegneria Matematica, Modelli Matematici e Simulazioni Numeriche, Torino, Italy.}
	\address{$^3$ Azienda Ospedaliera San Camillo, Unità Operativa Complessa di Cardiochirurgia e Chirurgia dei Trapianti, Roma, Italy.}
	\address{$^4$ Azienda Ospedaliera FBF Luigi Sacco, Dipartimento di Cardiochirurgia, Milano, Italy.}
	
	\email{mgirfogl@sissa.it}
	\email{lscandur@sissa.it}
	\email{fballarin@sissa.it}
	\email{giuseppe.infantino@studenti.polito.it}
	\email{nicolo\_francy84@hotmail.it}
	\email{andrea.montalto@libero.it}
	\email{grozza@sissa.it} 
	\email{roberto.scrofani@asst-fbf-sacco.it}         
	\email{marina.comisso@gmail.com}
	\email{fr.musumeci@gmail.com}
	
	\begin{abstract}
		Reduced order modeling (ROM) techniques are numerical methods that approximate the solution of parametric partial differential equation (PDE) by properly combining the high-fidelity solutions of the problem obtained for several configurations, i.e. for several properly chosen values of the physical/geometrical parameters characterizing the problem. 
		In this contribution, we propose an efficient non-intrusive data-driven framework involving ROM techniques in computational fluid dynamics (CFD) for hemodynamics applications. 
		By starting from a database of high-fidelity solutions related to a certain values of the parameters, we apply the proper orthogonal decomposition with interpolation (PODI) and then reconstruct the variables of interest for new values of the parameters, i.e. different values from the ones included in the database. Furthermore, we present a preliminary web application through which one can run the ROM with a very user-friendly approach, without the need of having expertise in the numerical analysis and scientific computing field. The case study we have chosen to test the efficiency of our algorithm is represented by the aortic blood flow pattern in presence of a Left Ventricular Assist Device (LVAD) when varying the pump flow rate.\\
		
	\end{abstract}
	
	\maketitle
	
	\section{Introduction}\label{sec:intro}
	Reduced order modeling (ROM) (see, e.g., \cite{ModelOrderReduction}) is a well-spread technique used both in academia and in industry. It has been introduced as an efficient tool to approximate full order systems by significantly reducing the computational cost required to obtain numerical solutions in a parametric setting. ROM consists in two main stages: an \emph{offline} phase that can be carried out on high performance computing facilities, and an \emph{online} one that hinges on a system of reduced dimensionality to perform the parametric computation on portable devices. In the \emph{offline} phase, the reduced order space is built starting from full order complex simulations computed for certain values of the physical and/or geometrical parameters. In this work, we employ the proper orthogonal decomposition (POD) for the detection of the reduced basis functions that span this new reduced space.
	After the creation of such a space, in the \emph{online} phase a new parametric solution is obtained as a linear combination of the precomputed reduced basis functions, by means of an interpolation carried out by using RBF functions \cite{Lazzaro2002}.
	The resulting ROM is thus called proper orthogonal decomposition with interpolation (PODI) \cite{Bui-Thanh2004}.
	
	The aim of this work is the development of an efficient non-intrusive data-driven reduced order model to be used within hemodynamics framework. The reader can find examples of the ROM application in the hemodynamics field in \cite{Ballarin2015,Ballarin2016,Tezzele2018,Zainib2020,Girfoglio2020}. We highlight that the online evaluation of the data-driven approach used here is based only on data and does not require knowledge about the governing equations that describe the system. It is also non-intrusive, i.e. no modification of the simulation software is carried out. For this reason it is particularly versatile thanks to its capability to be coupled with commercial solvers as well. It should be noted that many efforts are making in order to integrate ROM and technological innovation. From this viewpoint, a crucial step is the web server ARGOS \cite{Argos}, developed by mathLab group at SISSA that will make possible the exploitation of reduced order models to a wide category of people working in industrial and biomedical contexts. Through specific web applications the user will be able to solve many complex problems without the need of being an expert in numerical analysis and scientific computing. In particular, it is expected that the ATLAS project \cite{Atlas} will collect all cardiovascular applications. In this framework, we present a preliminary web application through which one can run the ROM by using a very user-friendly GUI interface. The benchmark we have chosen to test the efficiency of our algorithm is represented by the aortic blood flow pattern in presence of a Left Ventricular Assist Device (LVAD) (see, e.g., \cite{A,B,C,D,E,F,G}) when varying the pump flow rate (see, e.g., \cite{Girfoglio2020,Bazilevs2009,Mazzitelli2016}).
	
	The work is organized as follows. In Sec. \ref {sec:fom} we present the general parametric full order model governing hemodynamics problems, over which we apply the proposed numerical methodology. In Sec. \ref {sec:rom} we present the PODI method, whilst in Sec. \ref{sec:results} we show the numerical setting of the problem and the achieved results, as well as provide a brief description of the web application developed. Finally, in Sec. \ref{sec:conclusion} conclusions and perspectives are provided.
	
	\section{The full order model}\label{sec:fom}
	In this work we consider the blood as modeled by the unsteady incompressible Navier-Stokes equations described in
	an Eulerian framework. We consider a fixed domain $\Omega \subset \mathbb{R}^d$ with $d = 2, 3$ over a time
	interval of interest ($t_0$, $T$) $\subset \mathbb{R}^+$. Let $\bm{\pi} \in \mathcal{P} \subset \mathbb{R}^P$
	be a parameter vector in a $P$-dimensional parameter space $\mathcal{P}$. We have
	
	\begin{equation}\label{eq:NS1}
		\rho\, \partial_t \u(\bm{x},t; \bm{\pi})+ \rho\,\nabla \cdot \left(\u(\bm{x},t; \bm{\pi}) \otimes \u(\bm{x},t; \bm{\pi})\right) - 2\mu \Delta\u(\bm{x},t; \bm{\pi}) + \nabla p(\bm{x},t; \bm{\pi}) = 0,
	\end{equation}
	
	\begin{equation}\label{eq:NS2}
		\nabla \cdot \u(\bm{x},t; \bm{\pi}) = 0,
	\end{equation}\\
	in $\Omega \times [t_0,T]$, endowed with proper boundary conditions. In \eqref{eq:NS1}-\eqref{eq:NS2}, $\partial_t$ denotes the time derivative, $\rho = 1060$ kg/m$^3$ is the blood density, $\mu = 0.004$ Pa $\cdot$ s is the blood \emph{dynamic} viscosity, $\bm{u}$ is the blood velocity and $p$ is the pressure.
	
\begin{figure*}
		\centering
		\includegraphics[width=0.4\textwidth]{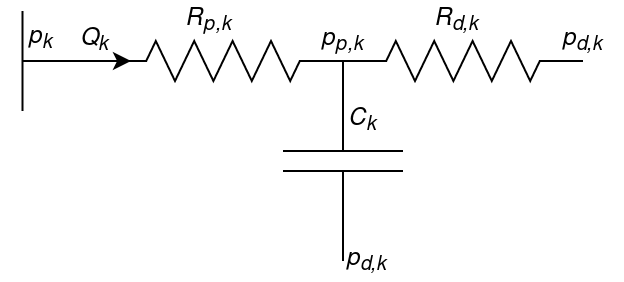}
		\caption{Three-element Windkessel model for the generic outlet $k$}
		\label{fig:RCR}
	\end{figure*}

	We impose a no slip boundary condition on the wall of the domain. At the inflow, we prescribe a known flow-rate and ${\partial p}/{\partial \n} = 0$ where $\n$ is the outward normal. On the other hand, in order to enforce realistic outflow boundary conditions at each outlet of the domain, we consider the Windkessel model based on the electric--hydraulic analogy \cite{Shi2011b}. By representing the blood pressure and flow rate with voltage and current, respectively, and by describing the effects of friction and inertia in blood flow and of vessel elasticity with resistance \emph{R}, inductance \emph{L} and capacitance \emph{C}, respectively, the methods for analysis of electric circuits can be borrowed and applied to the investigation of cardiovascular dynamics. In this work, we consider a three-element Windkessel RCR model \cite{Westerhof2008}. It consists of a proximal resistance $R_{p,k}$, a compliance $C_k$, and a distal resistance $R_{d,k}$, for each outlet $k$ (Fig. \ref{fig:RCR}).
	
	The downstream pressure, $p_k$, is expressed through the following differential algebraic equations (DAE) system:
	\begin{equation}\label{eq:RCR}
		\begin{cases}
			C_k \dfrac{dp_{p,k}}{dt} + \dfrac{p_{p,k} - p_{d,k}}{R_{d,k}} = Q_k, \\ \\
			p_k - p_{p,k} = R_{p,k}Q_k, \\
		\end{cases}
	\end{equation}\\
	where $Q_k$ is the flow rate, and $p_{p,k}$ and $p_{d,k}$ are the proximal and the distal pressure, respectively.
	
	For the space discretization of problems \eqref{eq:NS1}-\eqref{eq:NS2}, we adopt the Finite Volume (FV) approximation. A partitioned approach has been used to deal with the pressure-velocity coupling. In particular a Poisson equation for pressure has been used. This is obtained by taking the divergence of the momentum equation \eqref{eq:NS1} and exploiting the divergence free constraint \eqref{eq:NS2},
	
	\begin{equation}\label{eq:Poisson}
		\Delta p = -\nabla \cdot \left( \nabla \cdot \left(\u \otimes \u\right) \right).
	\end{equation}\\
	We have used the PISO algorithm \cite{PISO} employed in the finite volume C++ library OpenFOAM\textsuperscript{\textregistered} \cite{Weller1998}. For more details, we refer the reader to \cite{Girfoglio2020}.
	
	\section{The reduced order model}\label{sec:rom}
	The reduced order model we propose is the so-called \emph{proper orthogonal decomposition with interpolation}. It is a technique widely used within the reduced order modeling community in the study of parametric problems.
	POD allows to extract, from a set of high-dimensional snapshots, the optimal basis which minimizes the error between the original snapshots and their orthogonal projection. The data-driven approach used in this work is based only on data and does not require knowledge about the governing equations that describe the system (and which generated the snapshots). It is also non-intrusive, i.e., no modification of the simulation software is carried out. Still, there are works that use non-intrusive methods that are not data-driven (see, e.g., \cite{Zou2017}). The original snapshots are projected onto the POD space in order to reduce their dimensionality. Then the solution manifold is approximated using an interpolation technique. In this work, we will use a radial basis function (RBF) interpolation \cite{Lazzaro2002}. Several examples of applications based on this so-called POD with interpolation (PODI) \cite{Bui-Thanh2004} technique can be found in literature, in a wide range of contexts (see, e.g., \cite{Demo2019}).
	
    We are going to briefly describe the method that consists in two phases:\\
	\begin{itemize}
		\item[--] \emph{Offline}: let $N$ denote the number of degrees of freedom, e.g.\ associated to the FV discretization introduced in the previous section. Let $\bm{\varphi}_i$, with $i = 1, \dots , N_s$, be the snapshots related to a generic variable of interest collected by solving the high-fidelity problem, with different values of the input parameters $\bm{\pi}_i$, resulting in $N_s$ input-output pairs ($\bm{\pi}_i$, $\bm{\varphi}_i$). The snapshots matrix $\bm{S}$ is assembled by arranging the snapshots as columns, i.e.\ $\bm{S} = [\bm{\varphi}_1, \bm{\varphi}_2, \dots, \bm{\varphi}_{N_s}]$. By applying the singular value decomposition to this matrix, we have:
		
		\begin{equation}\label{eq:svd}
			\bm{S} = \bm{U}\bm{\Sigma}\bm{V}^* \approx \bm{U}_k\bm{\Sigma}_k{\bm{V}_k}^*,
		\end{equation}\\
		where $\bm{U} \in \mathcal{A}^{N \times N_s}$ is the unitary matrix containing the left-singular vectors, $\bm{\Sigma} \in \mathcal{A}^{N_s \times N_s}$  is the diagonal matrix containing the singular values $\lambda_i$, and $\bm{V} \in \mathcal{A}^{N_s \times N_s}$, with the symbol $^*$  denoting the conjugate transpose. The left-singular vectors, namely the columns of $\bm{U} = [\bm{\phi}_1, \bm{\phi}_2, \dots, \bm{\phi}_{N_s}]$, are the so-called POD modes.
        In order to reduce the dimensionality of the problem, we can keep the first $k$ modes to
		span the optimal space with dimension $k$ to represent the snapshots. By considering that
		the singular values are returned in decreasing order, we could truncate the
		number of modes simply selecting the first $k$ columns of $\bm{U}$.
		Therefore, the matrices $\bm{U}_k \in \mathcal{A}^{N \times k}$, $\bm{\Sigma}_k \in \mathcal{A}^{k \times k}$, $\bm{V}_k \in \mathcal{A}^{N_s \times k}$ in Eq. \ref{eq:svd} are the truncated matrices with rank $k$.
		
		After constructing the POD space, we can project the original snapshots onto this space. We compute $\bm{C} \in \mathcal{R}^{k \times N_s}$ as $\bm{C} = {\bm{U}_k}^T \bm{S}$, where
		the columns of $\bm{C}$ are the so-called modal coefficients. We express the input snapshots as a linear combination of the modes using such coefficients. Then, we have:
		
		\begin{equation}\label{eq:svd3}
			\bm{\varphi}_i \approx \sum_{j=1}^{k}  {\alpha}_j(\bm{\pi}_i)\bm{\phi}_j, \quad \forall i \in [1,2, \dots, N_s],
		\end{equation}\\
		where ${\alpha}_j(\bm{\pi}_i)$ are the elements of $\bm{C}$.
		At a given mode $\bm{\phi}_j$, the ($\bm{\pi}_i$, ${\alpha}_j(\bm{\pi}_i)$) pairs sample the solution manifold in the parametric space. The interpolation of the modal coefficients $\alpha_j(\bm{\pi}_i)$ in the parameter space is carried by using RBF functions. It is based on the following formula:
			
			\begin{equation}\label{eq:online1}
				A_j(\bm{\pi}) = \sum_{m=1}^{N_s} w_{j,m} \zeta_{j,m} \left(||\bm{\pi} - \bm{\pi}_m||_{L^2(\mathbb{R}^{\mathcal{P}+1})}\right),
			\end{equation}\\
			where $w_{j,m}$ are proper weights and $\zeta_{j,m}$ are the RBF functions which are chosen to be Gaussian functions, centered in $\bm{\pi}_m$.
		
		For the computation of the weights $w_{j,m}$, the following property has to be used:
			
			\begin{equation}\label{eq:online2}
				A_j(\bm{\pi}_i) = {\alpha}_j(\bm{\pi}_i),
			\end{equation}\\
			The last equation can be rewritten in form of a linear system:
			
			\begin{equation}\label{eq:online3}
				\bm{A}_j^{\zeta} \bm{w}_j = \bm{\alpha}_j.
			\end{equation}\\
			Thus, one could solve the latter linear system to obtain the weights $\bm{w}_j$ related to the mode $\bm{\phi}_j$ which will be stored to be then used in the online stage.\\
		\item[--]\emph{Online}: we are able, for any new parameter value $\bm{\pi^\star}$ to calculate the new coefficients $\alpha_j(\bm{\pi^\star})$, which are given simply by:
			
			\begin{equation}\label{eq:online4}
				\alpha_j(\bm{\pi^\star}) = \sum_{m=1}^{N_s} w_{j,m} \zeta_{j,m} \left(||\bm{\pi^\star} - \bm{\pi}_m||_{L^2(\mathbb{R}^{\mathcal{P}+1})}\right),
		\end{equation}\\
		Then, we compute the high-dimensional solution by projecting back the (approximated) modal coefficients to the original space:
		
		\begin{equation}\label{eq:online5}
			\bm{\varphi}(\bm{\pi^\star}) = \sum_{j=1}^{k}  \alpha_j(\bm{\pi^\star}) \phi_j
		\end{equation}\\
	\end{itemize}
	
	We remark that the procedure can be repeated for several variables of interests. Furthermore, it is not necessary for such a variable to be an unknown of the original system (such as velocity and pressure); indeed, we will use the PODI technique not only for primal quantities, but also for derived quantities such as wall shear stress (WSS).
	
	Regarding the technical implementation of the PODI method, we use the Python package called EZyRB \cite{eazyrb}.

	\section{Numerical results and discussion}\label{sec:results}
	In order to test the performance of the presented computational pipeline, we investigate the aortic blood flow pattern in presence of a Left Ventricular Assist Device (LVAD) when varying the pump flow rate. This case study has been thoroughly discussed both at high-fidelity (or full order model, FOM) and ROM level in \cite{Girfoglio2020}. Here, after summarizing some relevant computational details related to the clinical data, geometrical model and full order simulations (Sec. \ref{sec:results0}), we are going to further extend the ROM investigation with additional tests (Sec. \ref{sec:results1}) and to provide a brief description of the preliminary web application under development (Sec. \ref{sec:results2}).
	
	\subsection{Computational details}\label{sec:results0}
\begin{figure*}[]
		\centering
		\begin{overpic}[width=0.4\textwidth]{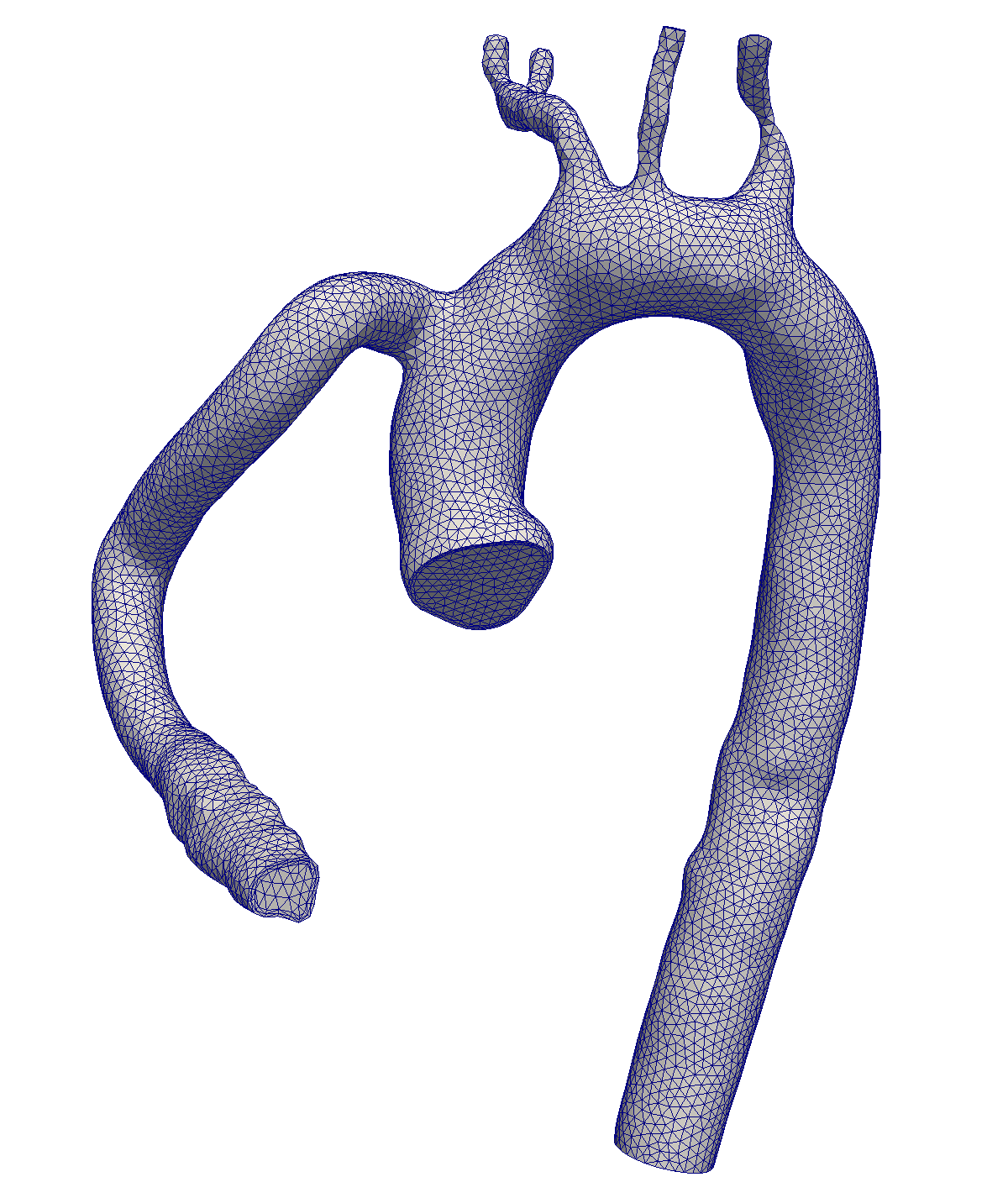}
			\put(30,100){\small{a)}}
		\end{overpic}
		\begin{overpic}[width=0.4\textwidth]{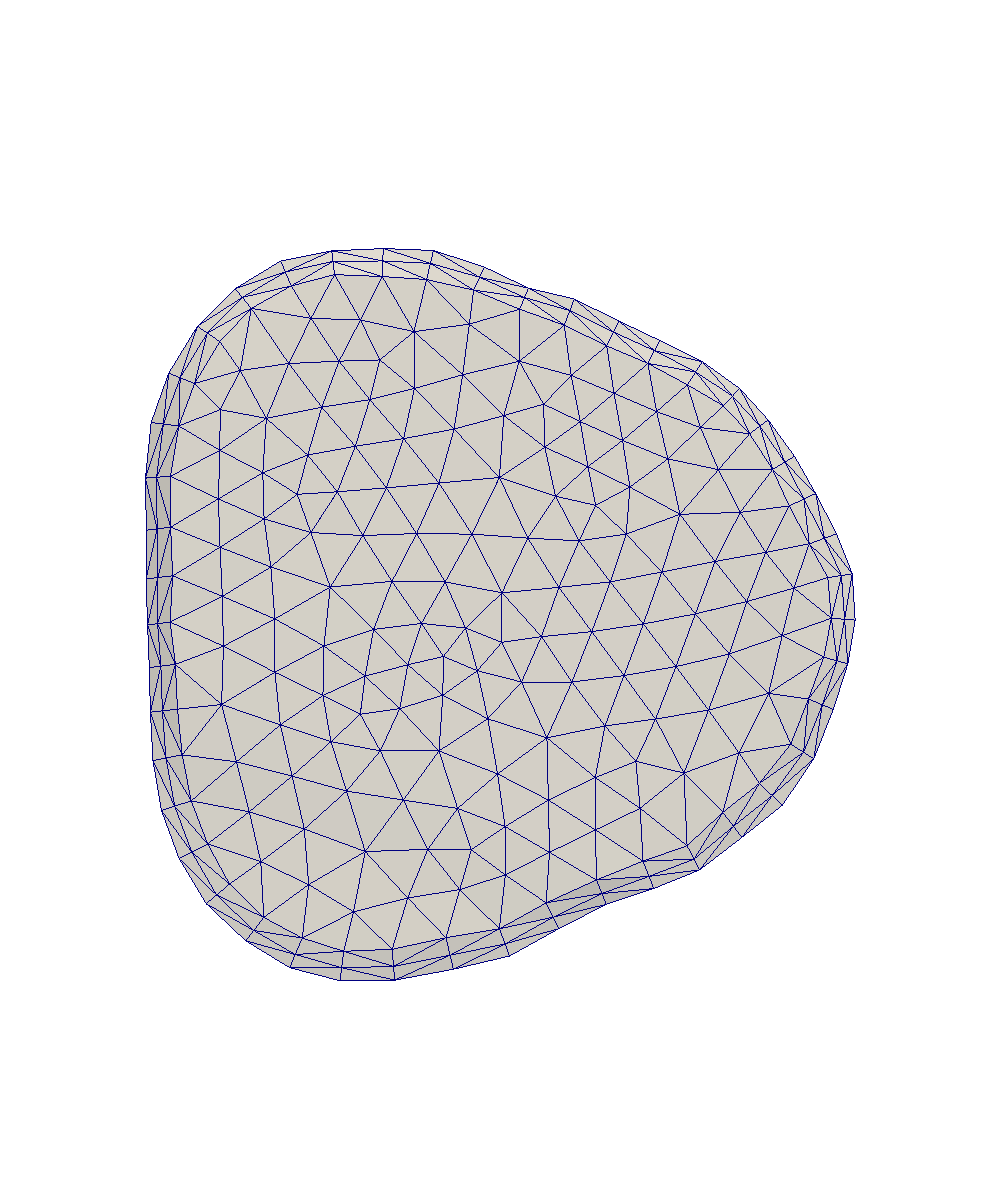}
			\put(30,95){\small{b)}}
		\end{overpic}
		\caption{View of the mesh: a) overview, b) a section next to the aortic inlet.}\label{fig:mesh}
	\end{figure*}

	A real patient-specific aorta model was reconstructed from Computed Tomography (CT) images by using the open source medical image analysis software 3D Slicer\textsuperscript{\textregistered} (\url{http://www.slicer.org}). The model is referred to the post-surgery configuration, i.e. after the implantation of the LVAD device (the Heartmate 3$^{TM}$ Left Ventricular Assist System \cite{HeartMate}), and includes the outflow cannula of the LVAD device (with inlet in the bottom left of Fig. \ref{fig:mesh}a), the ascending aorta (with inlet in the center of Fig. \ref{fig:mesh}a, which is observed from below in Fig. \ref{fig:mesh}b), brachiocephalic artery, right subclavian artery, right common carotid artery, left common carotid artery, left subclavian artery and descending aorta, as shown in Fig. \ref{fig:mesh}. We consider a tetrahedral computational grid with $h_{min} = 5.83e-4$ and $h_{max} = 3.2e-3$ for a total of $200k$ cells. The quality of this mesh is suitable for a FV solver: it features very low values of average non-orthogonality (30$^\circ$) and skewness (around 1). Fig. \ref{fig:mesh} shows the mesh. It should be noted that such a computational grid has been used in \cite{Girfoglio2020} where a mesh convergence analysis is carried out.
	
	By making reference to the equations \eqref{eq:NS1}-\eqref{eq:Poisson}, the convective term is discretized by using a first order upwind scheme. On the other hand, for the diffusive term,  a central differencing interpolation scheme with non-orthogonality correction is preferred. Regarding the pressure gradient, we use a linear interpolation scheme. For more details about such schemes, the reader can make reference to \cite{Weller1998,jasak1996error}. Finally, to discretize in time the equations \eqref{eq:NS1} and \eqref{eq:RCR}, we adopt Backward Differentiation Formula of order 1 (BDF1),  see e.g. \cite{quarteroni2007numerical}.
	
	Coefficients values of the Windkessel models used for the enforcement of the outlet boundary conditions shown in Table \ref{tab:bc_post} are based on \cite{Girfoglio2020}.
	
	Although the high-order simulations are obtained through time-stepping by solving the governing equations \eqref{eq:NS1}-\eqref{eq:Poisson}, we are interested in computing and collecting the steady states, i.e. solutions where $\partial_t \u$ vanishes.
	
	\begin{table*}
		\begin{center}
			\begin{tabular}{|c|c|c|c|}
				\hline
				& $R_{p,k}$ [dyne $\cdot$ s/cm$^5$] & $R_{d,k}$ [dyne $\cdot$ s/cm$^5$] & $C_k$ [cm$^5$/dyne]
				\\\hline
				Right subclavian artery      &  2.56e3 & 4.32e4  &  3.26e-5    \\\cline{1-4}
				Right common carotid artery       &  1.63e3 & 2.74e4 & 5.16e-5       \\\cline{1-4}
				Left common carotid artery        &  2.38e3 & 4e4 & 3.52e-5      \\\cline{1-4}
				Left subclavian artery   &  8.96e2 & 1.51e4  & 9.35e-5     \\\cline{1-4}
				Descending aorta     & 1.08e2  & 1.83e3 & 7.72e-4 \\\hline
			\end{tabular}
		\end{center}
		\caption{Windkessel coefficients: proximal resistance $R_{p,k}$ and distal resistance $R_{d,k}$, and compliance $C_k$, for each outlet $k$ \cite{Girfoglio2020}.}\label{tab:bc_post}
	\end{table*}
	
	\subsection{ROM results}\label{sec:results1}
	To train the ROM, we consider the range of pump flow rate (i.e., the flow rate at the inlet of outflow cannula) $PF \in [3, 5]$ that covers typical clinical values. Thus, we consider as parameter $\bm{\pi}$ the pump flow rate $PF$. In particular, we choose equispaced distributions inside the ranges $PF \in [3, 3.8]$ and $PF \in [4.2, 5]$. The sampling frequency is 0.2 for both ranges, so that we have a database including 10 snapshots related to the high-fidelity steady state solutions. It should be noted that in \cite{Girfoglio2020} we have showed that, for this benchmark, the number of snapshots does not affect significantly the accuracy of the ROM.
	One new value of $PF$ in which the ROM has not been trained but which belongs to the range of the training space, $PF = 4$, is used to evaluate the performance of the parametrized ROM. In \cite{Girfoglio2020} a equispaced distribution of 11 snapshots inside the range $PF \in [3, 5]$ was used, i.e. the snapshot related $PF = 4$ was included in the FOM database, and ROM was performed for $PF = 3.45$ and $PF = 4.35$. We note that, in \cite{Girfoglio2020}, the distance between the parameter values for which the ROM is performed and the nearest snapshot is 0.05. On the other hand, here such a distance is larger, 0.2.
	
	Table \ref{tab:cum_eigen} shows the cumulative energy of the eigenvalues for pressure $p$, wall shear stress WSS, and velocity components, $u_x$, $u_y$ and $u_z$.
	\begin{table*}
		\begin{center}
			\begin{tabular}{|c|c|c|c|c|c|}
				\hline
				N & $p$ & WSS &  $u_x$ & $u_y$ & $u_z$ \\\hline
				1 &  0.9999  & 0.9899  &  0.9834 & 0.9729 & 0.9785\\\hline
				2 &  0.9999  &  0.9957 &  0.9949 & 0.9903 & 0.9923 \\\hline
			\end{tabular}
		\end{center}
		\caption{Cumulative energy of the eigenvalues for pressure $p$, wall shear stress WSS, and velocity components, $u_x$, $u_y$ and $u_z$.}\label{tab:cum_eigen}
	\end{table*}
	In order to retain the 99\% of the system's energy, 1 mode for $p$, 1 for WSS, 2 for $u_x$,  $u_y$ and $u_z$ are selected. It has been verified that considering a larger number of POD modes does not increase the accuracy of the ROM. To provide some quantitative results, the relative error in the $L^2$-norm, calculated as
	
	\begin{equation}\label{eq:error1}
		E_X = 100\dfrac{||X_{FOM} - X_{ROM}||_{L^2(\Omega)}}{||X_{FOM}||_{L^2(\Omega)}} \%
	\end{equation}\\
	where $X_{FOM}$ is the value of a particular field in the FOM model, and $X_{ROM}$ the one that is calculated using the ROM, is considered. In Table \ref{tab:relative_errors}, the relative error for all the variables of interest is reported.
	
	\begin{table*}
		\begin{center}
			\begin{tabular}{|c|c|c|c|c|c|}
				\hline
				& $p$ & WSS &  $u_x$ & $u_y$ & $u_z$ \\\hline
				$E_X$ &  0.5\% & 7.7\% & 8.5\% & 12.2\% & 11.4\% \\\hline
			\end{tabular}
		\end{center}
		\caption{$L^2$ norm relative errors for pressure $p$, wall shear stress WSS, and velocity components, $u_x$, $u_y$ and $u_z$, for $PF = 4$ l/min.}\label{tab:relative_errors}
	\end{table*}
	
\begin{figure*}[h]
		\centering
		\begin{overpic}[width=0.3\textwidth]{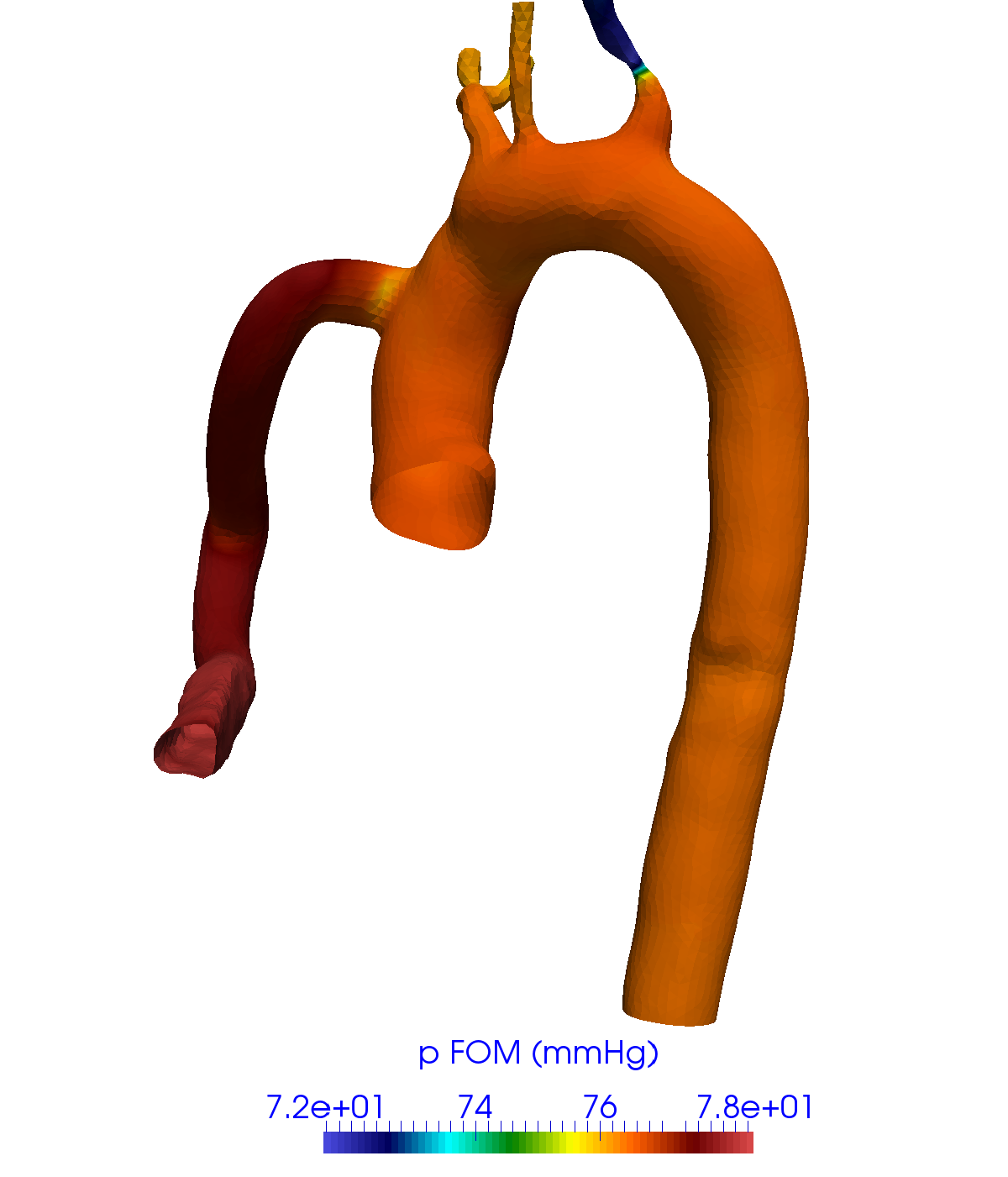}
		\end{overpic}
		\begin{overpic}[width=0.32\textwidth]{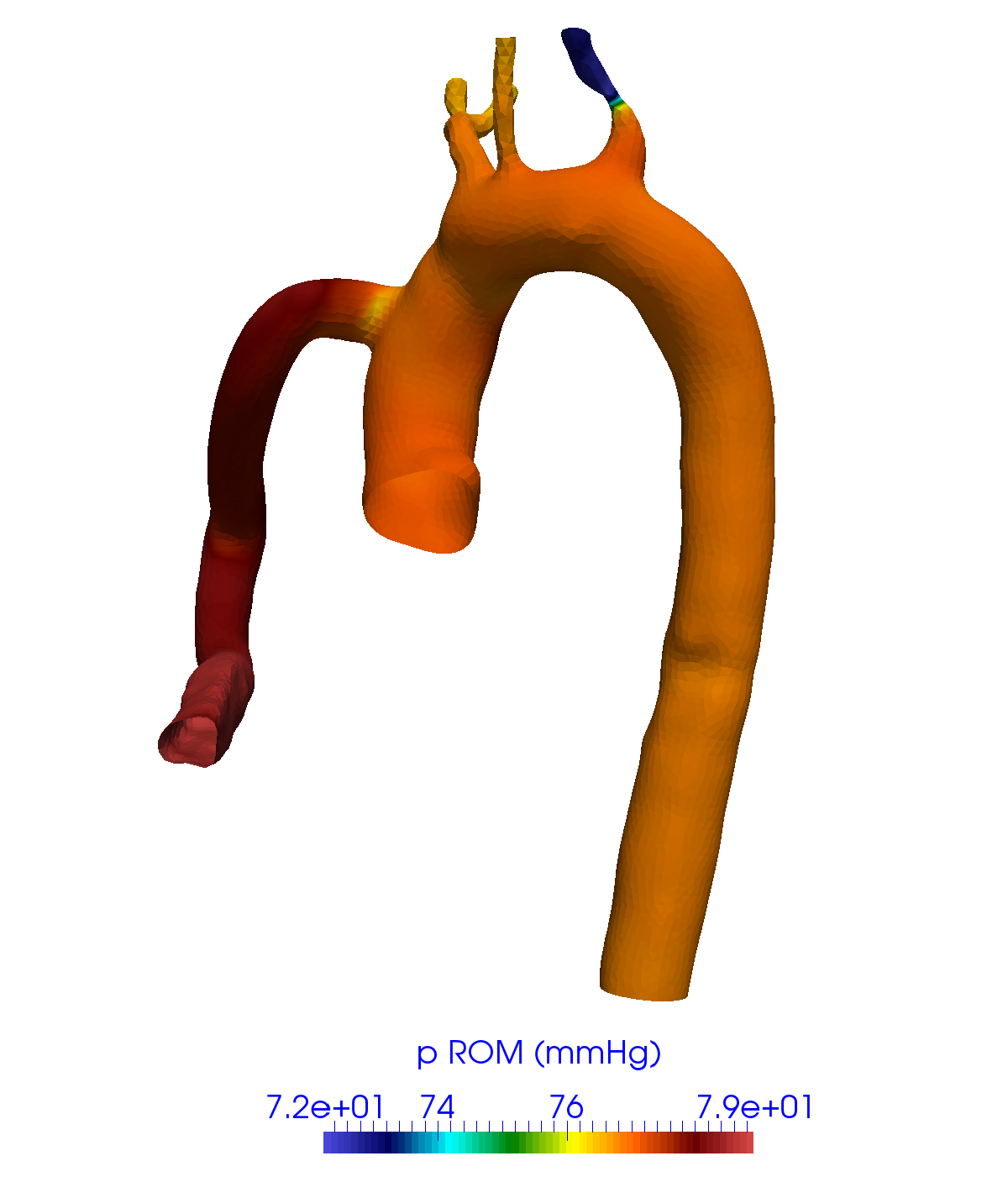}
		\end{overpic}
		\begin{overpic}[width=0.33\textwidth]{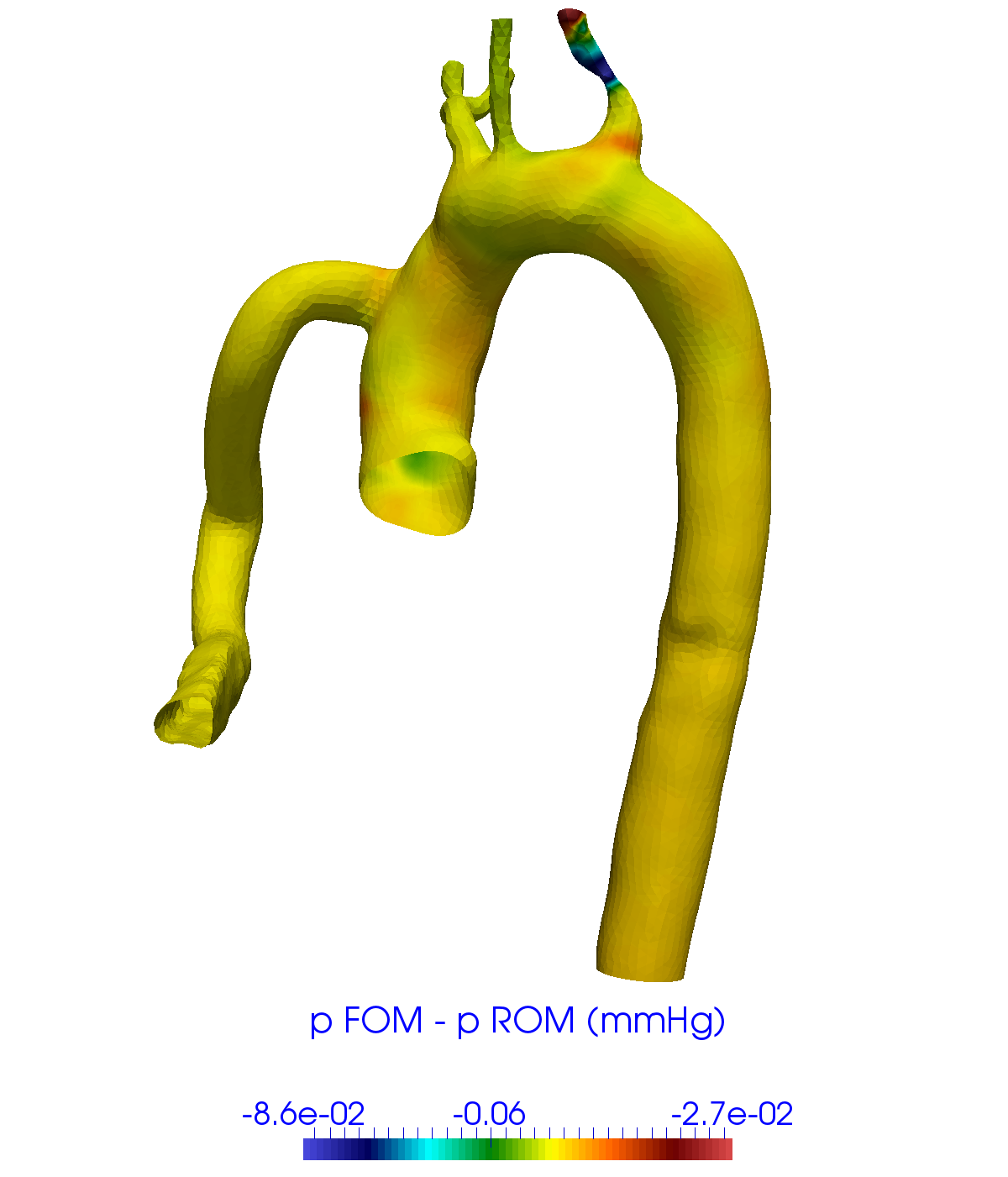}
		\end{overpic}
		\begin{overpic}[width=0.3\textwidth]{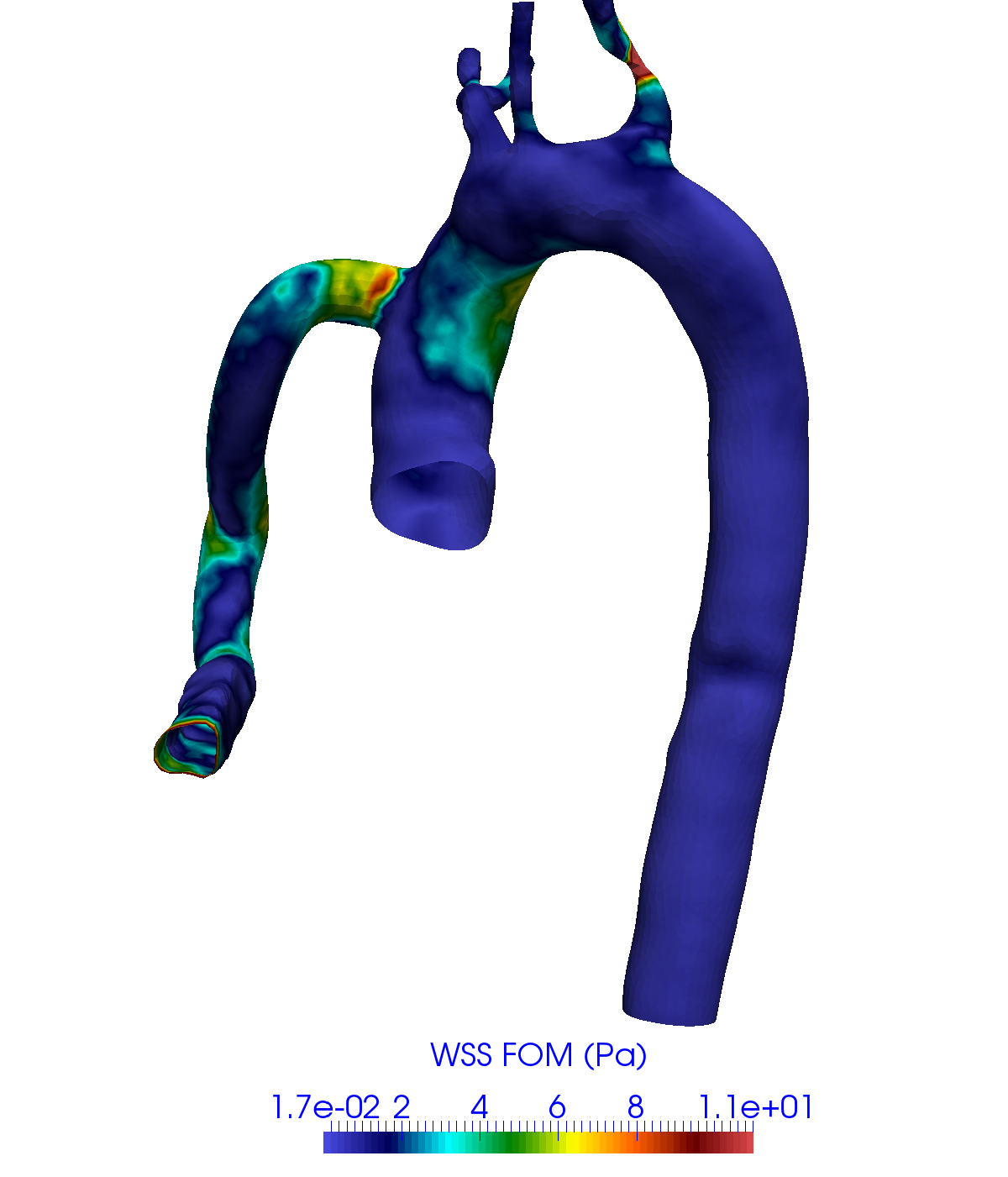}
		\end{overpic}
		\begin{overpic}[width=0.3\textwidth]{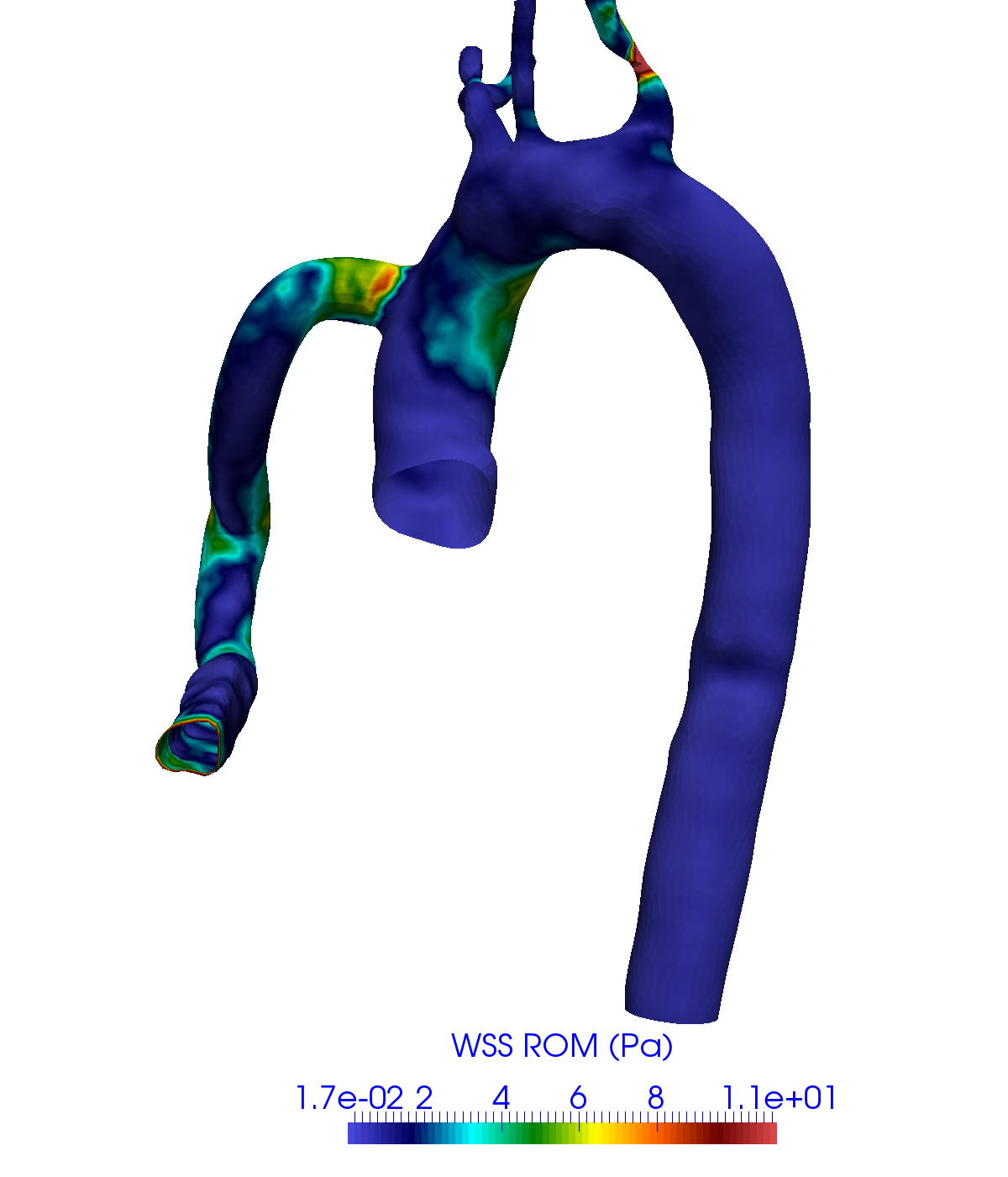}
		\end{overpic}
		\begin{overpic}[width=0.3\textwidth]{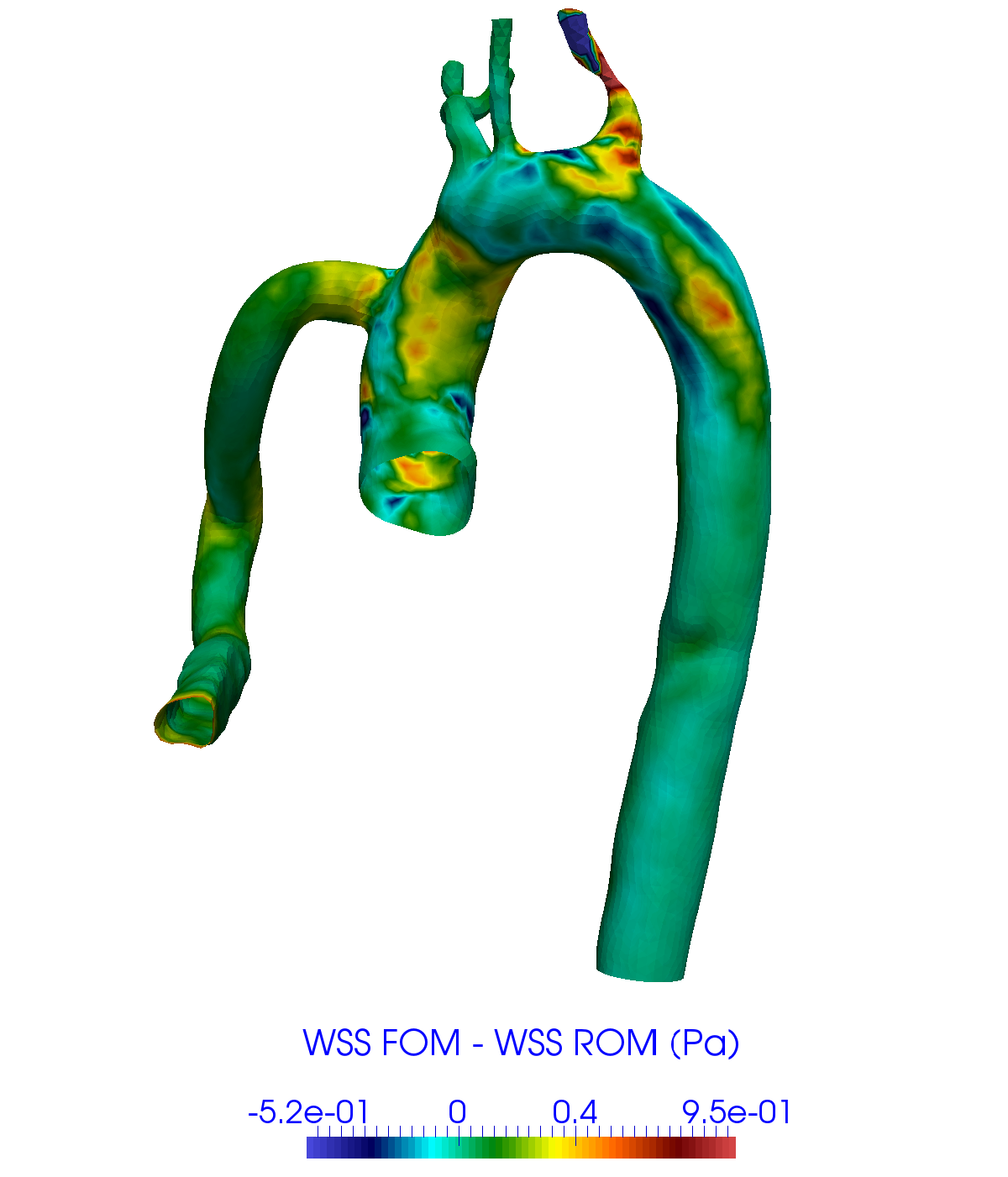}
		\end{overpic}
		\begin{overpic}[width=0.32\textwidth]{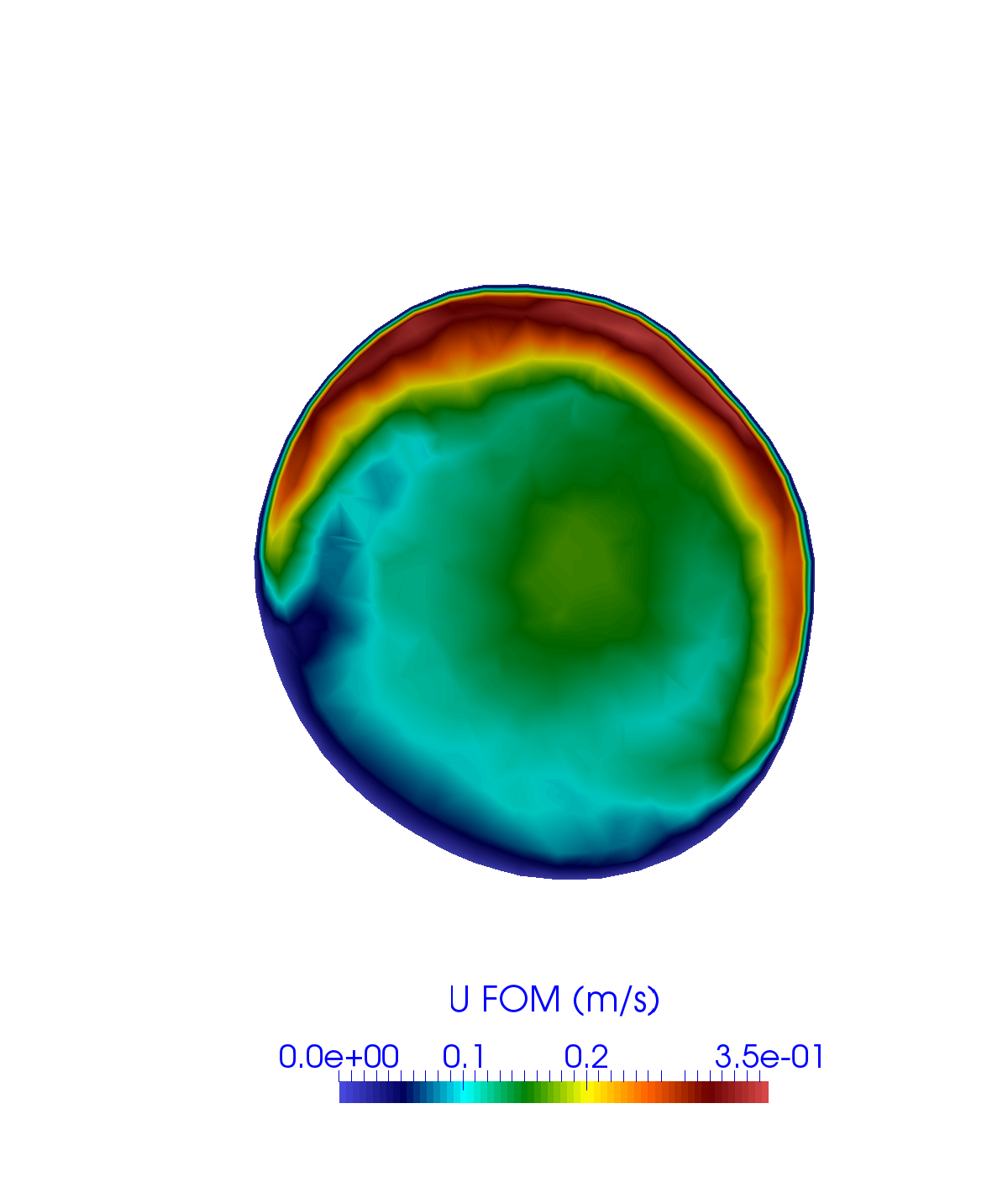}
		\end{overpic}
		\begin{overpic}[width=0.29\textwidth]{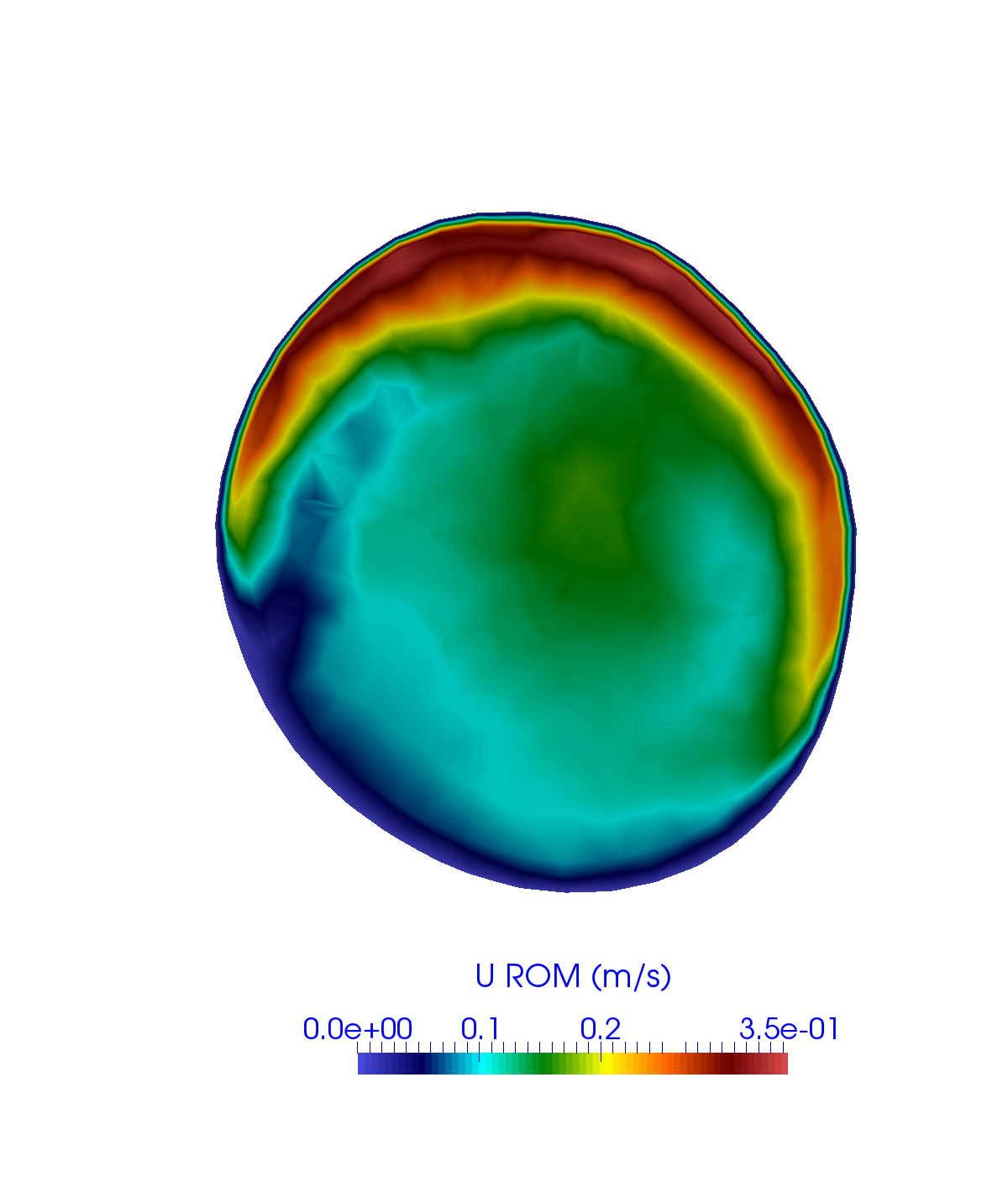}
		\end{overpic}
		\begin{overpic}[width=0.30\textwidth]{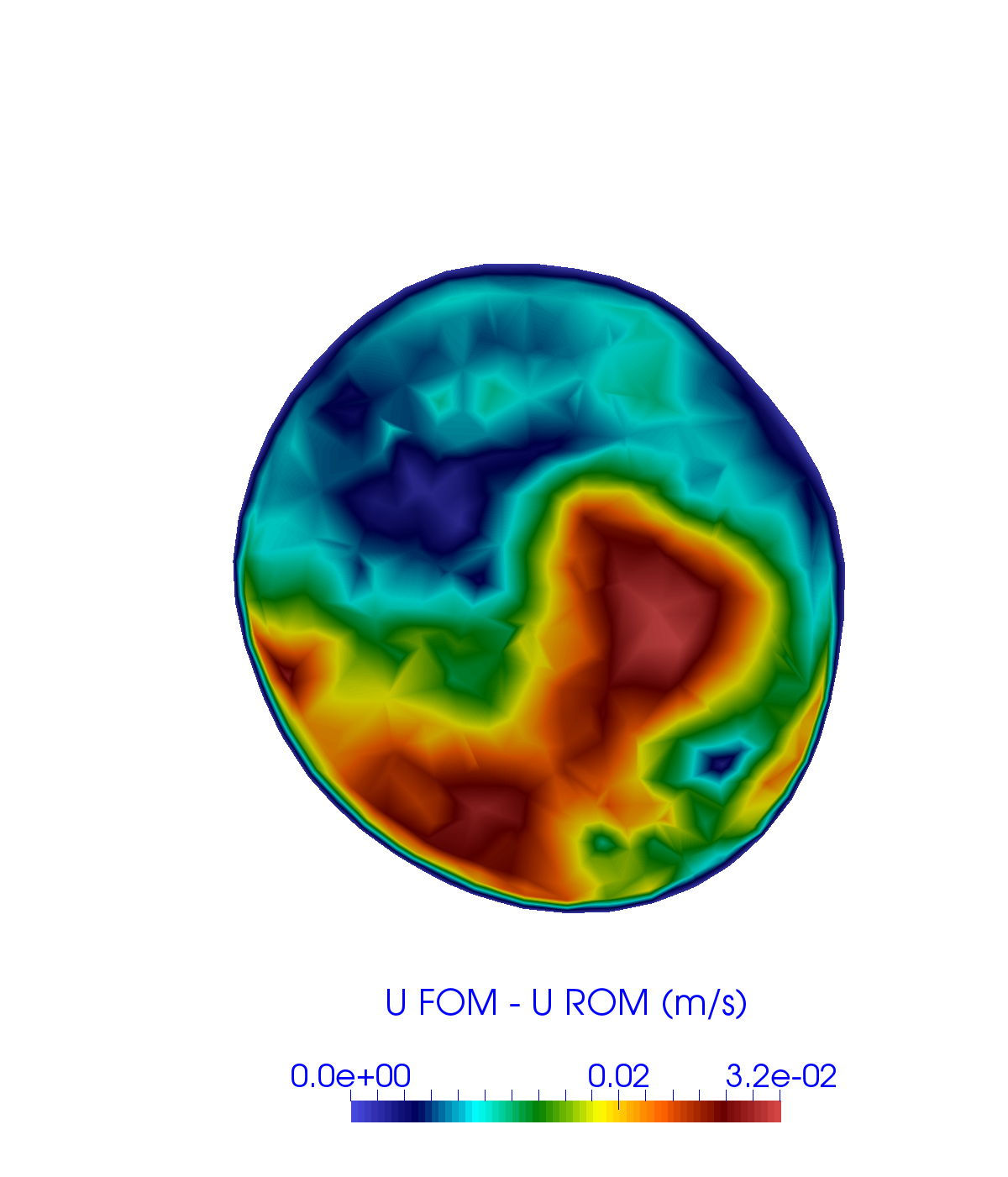}
		\end{overpic}
		\caption{Comparison of the FOM/ROM pressure (1st row), WSS (2nd row) and velocity steady-state solutions related to a section of the ascending aorta next to the anastomosis location (3nd row) at $PF = 4$.}\label{fig:ROM_4}
	\end{figure*}

	Fig. \ref{fig:ROM_4} displays a comparison between FOM and ROM for $p$ and WSS, and for the velocity related to a section of the ascending aorta next to the anastomosis location. The comparison indicates that the ROM is able to provide a good reconstruction for all the variables.
	
	Finally, we comment on the computational costs. The CPU time required by a FOM simulation is 9600s and the one of the ROM, that is related to the computation of the modal coefficients and reconstruction of the fields, is 40s. This corresponds to a speed-up of $\approx 240$, that demonstrates the fact that it is possible to use the ROM in the place of the FOM in order to obtain accurate simulations with a significant reduction of the computational cost.
	
	\subsection{Web application}\label{sec:results2}
	Due to the aforementioned speedup, research activities based on techniques (e.g., ROMs) that leads to technological innovation for real time calculation is acquiring considerable relevance and popularity also in the biomedical field.
	The combination of ROMs with technological development through a web interface would allow real time data to be accessed in hospitals and operating rooms on portable devices. In this scenario, the web server ARGOS \cite{Argos} has been created, which has the task of proposing a platform to favor a more widespread exploitation of real time computing through a simple ``click''. It is a very intuitive and smooth web platform, which does not need a strong experience in numerical analysis, fluid dynamics or scientific computing field to be used. ARGOS offers a wide variety of applications related to several problems and in particular it contains the section ATLAS \cite{Atlas} focused on the cardiovascular field. Here, we are going to provide a brief description of the web application under development which has the aim to support the user (which in this case could either be a scientist involved in the manufacturing of the pump, or a medical doctor interested in evaluation the hemodynamics in different operating scenarios) to set the LVAD device according to the need of the patient.
	
\begin{figure*}[]
		\centering
		\begin{overpic}[width=.99\textwidth]{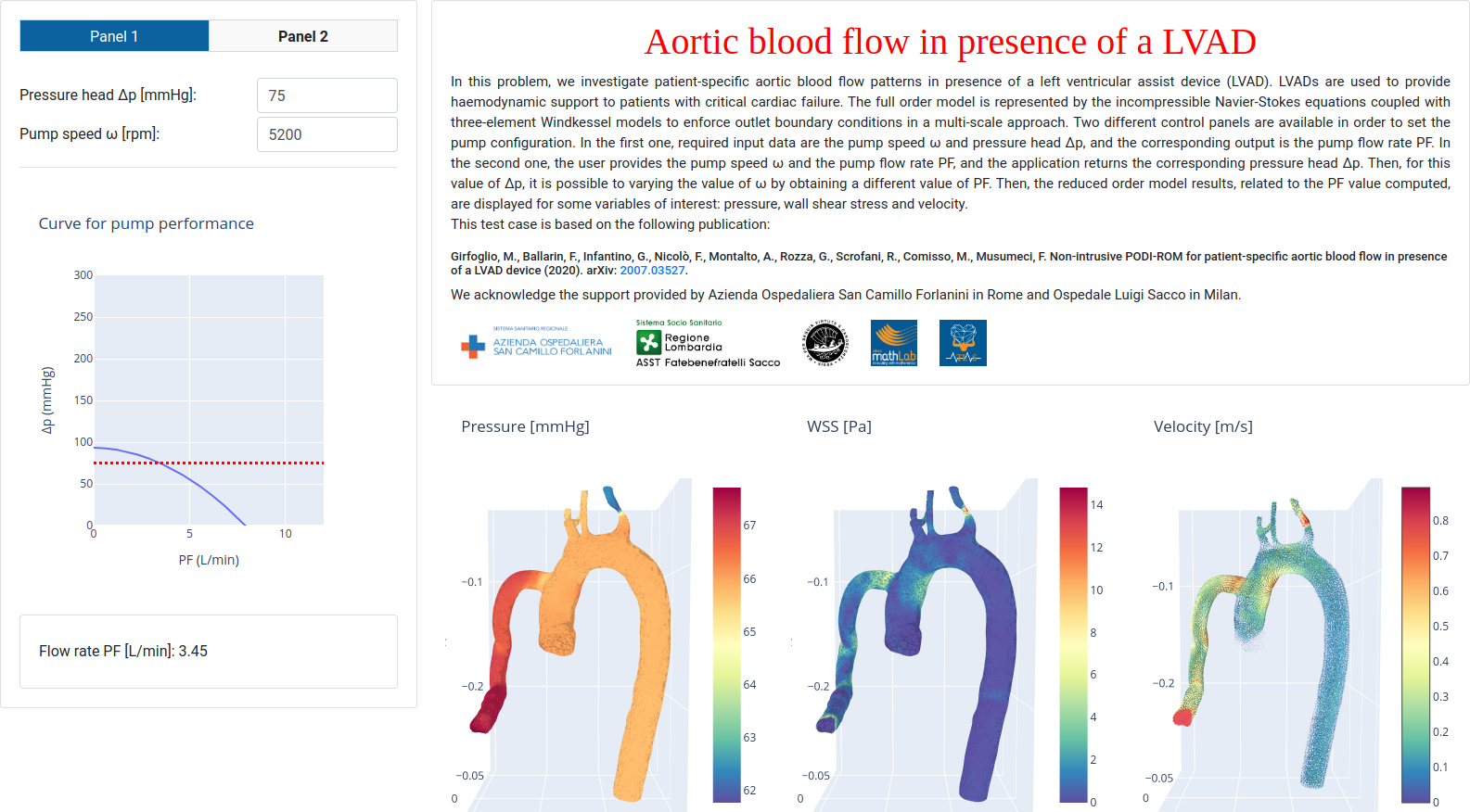}
			\put(-3,52){\small{a)}}
		\end{overpic}
		\begin{overpic}[width=.40\textwidth]{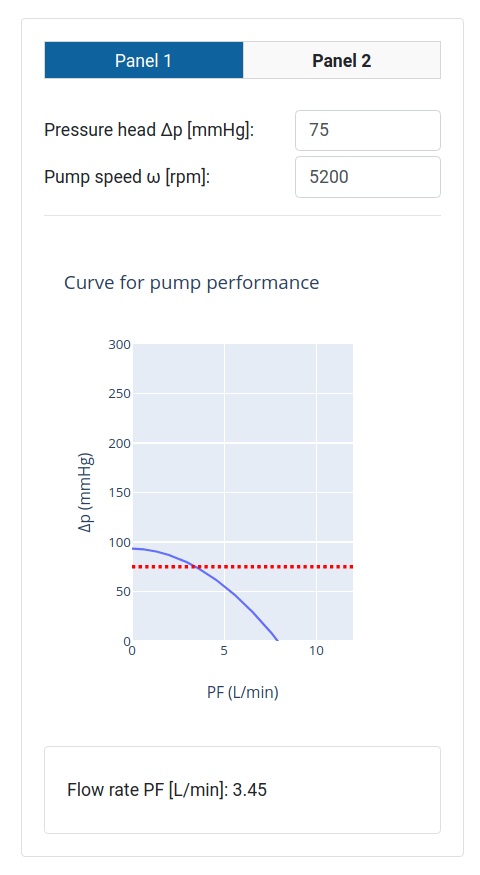}
			\put(-2,93){\small{b)}}
		\end{overpic}\hspace{1cm}
		\begin{overpic}[width=.341\textwidth]{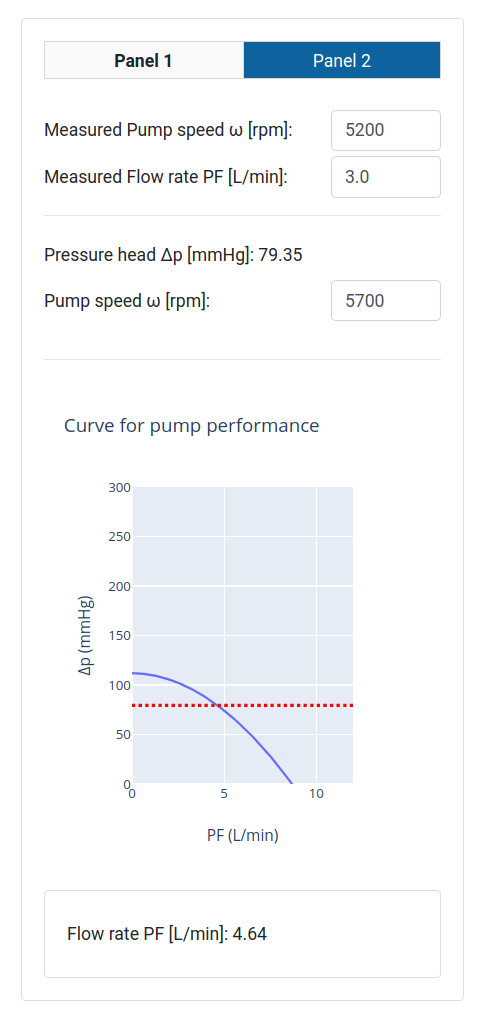}
			\put(-2,94.1){\small{c)}}
		\end{overpic}
		\caption{LVAD web application: a) overview, b) close-up of Panel 1, c) close-up of Panel 2.}\label{fig:app}
	\end{figure*}

\begin{figure*}
		\centering
		\includegraphics[width=0.5\textwidth]{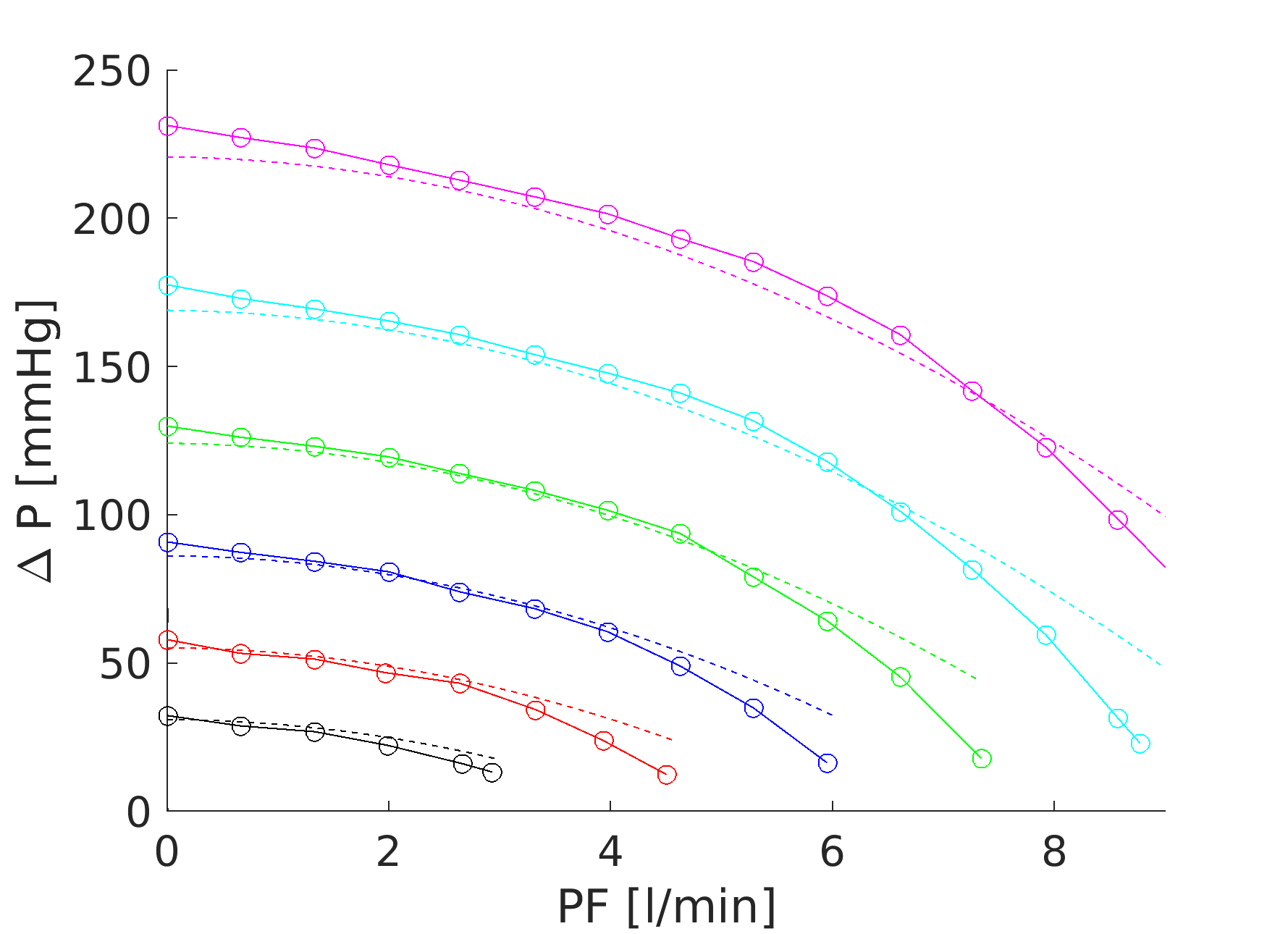}
		\caption{Pressure head ($\Delta P$) - volume flow rate ($PF$) curves (continuous line with circles) and analytical fitting (dashed line) based on eq. \ref{eq:pump} for Heartmate 3$^\text{TM}$ \cite{HeartMate} pump at several pump speed values: $\omega = 3000$  rpm (black), $\omega = 4000$ rpm (red), $\omega = 5000$ rpm (blue), $\omega = 6000$ rpm (green), $\omega = 7000$ rpm (cyan), and $\omega = 8000$ rpm (magenta).}
		\label{fig:pump_dynamics}
	\end{figure*}

	Figure \ref{fig:app} a) displays a screenshot of the application under development \cite{Lvad}.
	On the left side, the user can set up the pump configuration. We have two control panels related to two different settings denoted as Panel 1 (Figure \ref{fig:app} b)) and Panel 2 (Figure \ref{fig:app} c)). In Panel 1, required input data are the pressure head $\Delta P$ and pump speed $\omega$, and the corresponding output is the pump flow rate $PF$. The target user for this panel is a scientist involved in the design of the pump (i.e., for instance, by changing the pressure head $\Delta P$). In Panel 2, the user provides a measured pump speed $\omega$ and the corresponding measured pump flow rate $PF$, and the application returns the corresponding pressure head $\Delta P$. Then, for this value of $\Delta P$, it is possible to vary the value of $\omega$ to obtain a different value of $PF$. The target user for this panel is a medical doctor, who reads measured values of pump speed $\omega$ and flow rate $PF$ during an LVAD ramp test, and is then interested in predicting possible hemodynamics outcomes when changing the pump configuration set during the ramp test. The relationship between $\omega$, $\Delta P$ and $PF$ used in the application is given by the following analytical relationship
	
	\begin{align}\label{eq:pump}
		\Delta P = K_A\cdot\omega^2+K_B\cdot\omega\cdot PF+K_C\cdot PF^2
	\end{align}\\
	where $K_A$, $K_B$ and $K_C$ are constants which depend on pump design (Table \ref{tab:coeffs_fit}), that provides an acceptable fit as showed in Fig. \ref{fig:pump_dynamics}.
	The $\Delta P-PF$ analytical curve also is displayed in the application, below the control panels, and updates in real time when the pump speed $\omega$ is changed by the user. Currently, the application can be used only for $PF \in [3,5]$ because the FOM snapshots are collected for such a range of values. Then, if the values of $\omega$ and $\Delta P$ selected are such that the corresponding $PF$ is outside the range [3,5], an error message will appear. On the right part of the screen, we report a brief description of the web application, including relevant references, logos and acknowledgements, and we visualize in real time the solutions provided by the ROM related to the $PF$ value corresponding to the pump setting for the variables of interest: pressure, velocity and wall shear stress.
	
	\begin{table*}
		\begin{center}
			\begin{tabular}{|c|c|c|}
				\hline
				$K_A$ [mmHg/rpm$^2$] & $K_B$ [mmHg $\cdot$ l/min/rpm] & $K_C$ [mmHg $\cdot$ l$^2$/rpm$^2$] \\
				\hline
				3.45e-6      & -5.9e-5  &  -1.45 \\\hline
			\end{tabular}
		\end{center}
		\caption{Constants of the analytical fitting for the pump dynamics (see equation (\ref{eq:pump})).}\label{tab:coeffs_fit}
	\end{table*}
	
	\section{Conclusion and perspectives}\label{sec:conclusion}
	In this work, an efficient non-intrusive data-driven reduced order modelling to be used within hemodynamics framework is presented. The FOM is represented by the incompressible Navier-Stokes equations discretized by using a FV technique. Furthermore, the development of the ROM is carried out by using the proper orthogonal decomposition with interpolation (PODI). The online phase of the ROM results in a data-driven approach which is based only on data and does not require knowledge about the governing equations that describe the system. It is also non-intrusive, i.e. no modification of the simulation software is required. For this reason it is particularly versatile thanks to its capability to be coupled also with commercial solvers. Moreover, we have presented a preliminary web application through which one can run the ROM by using a very user-friendly interface, without the need of having a specific numerical expertise, and thus possibly widening the use of numerical tools to practitioners. The benchmark we have chosen to test the efficiency of our algorithm is represented by the aortic blood flow pattern in presence of a Left Ventricular Assist Device (LVAD) when varying the pump flow rate. We show that the ROM provides accurate solutions with a significant reduction of the computational cost, up to at least two orders of magnitudes.
	
	As a follow-up of the present work, we are going to make further efforts in order to improve the ease of use of the web application. We are also moving towards geometrical parametrization in the context of patient-specific geometries, extending e.g. the work carried out in \cite{Ballarin2016} to different problems and different model reduction techniques. Finally, we are interested in improving the full order model, by considering turbulence effects (see, e.g., \cite{Girfoglio2019,Girfoglio2020_b}), as well as coupling the fluid model with an elasticity model to simulate fluid-structure interaction (FSI) (see, e.g., \cite{Girfoglio2021,Bazilevs2009,Li2019,Mu2019}). This would make the presented pipeline more complete and versatile.
	
	\section*{Acknowledgements}\label{sec:acknowledgements}
	We acknowledge the collaboration with Mr Nicola Demo (SISSA mathLab) in the development of the ARGOS web server \cite{Argos}.
	
	We acknowledge the support provided by the European Research Council
	Executive Agency by the Consolidator Grant project AROMA-CFD ``Advanced
	Reduced Order Methods with Applications in Computational Fluid Dynamics'' -
	GA 681447, H2020-ERC CoG 2015 AROMA-CFD and INdAM-GNCS 2020 project ``Tecniche Numeriche Avanzate per Applicazioni Industriali''.
	
	\bibliographystyle{abbrv}
	\bibliography{bib/bibfile}
	
\end{document}